\documentclass[11pt]{amsart}

\usepackage{amsmath,amssymb,amsthm} 
\usepackage[mathscr]{eucal} 
\usepackage{latexsym}
\theoremstyle{plain}
\usepackage{enumerate}
\usepackage{hyperref}
\usepackage{color}

\usepackage[T1]{fontenc}
\usepackage[bitstream-charter]{mathdesign}
\usepackage{eucal}

\newtheorem{thm}{Theorem}
\newtheorem{lem}[thm]{Lemma}
\newtheorem{rem}{Remark}

\newcommand{\norm}[1]{{\left\Vert #1 \right\Vert}}
\newcommand{\set}[1]{{\left\lbrace #1 \right\rbrace}}

\DeclareMathOperator*{\divg}{div}

\title{
An improved Liouville-type theorem for the stationary tropical climate model
}

\author{Youseung Cho}

\author{Hyunjin In}

\author{Minsuk Yang}

\address{Yonsei University, 
Department of Mathematics, Yonseiro 50,
Seodaemungu, Seoul 03722, Republic of Korea}

\email{m.yang@yonsei.ac.kr}

\begin{document} 

\begin{abstract} 
In this paper, we study the Liouville-type property for smooth solutions to the steady 3D tropical climate model. 
We prove that if a smooth solution $(u,v,\theta)$ satisfies $u \in L^3 (\mathbb{R}^3)$, $v \in L^2 (\mathbb{R}^3)$, and $\nabla \theta \in L^2 (\mathbb{R}^3)$, then $u=v=0$ and $\theta$ is constant, which improves the previous result, Theorem 1.3 (Math. Methods Appl. Sci. 44, 2021)
by Ding and Wu.
\end{abstract}

\thanks{ 
\endgraf
}
\thanks{This work was supported by the National Research Foundation of Korea(NRF) grant funded by the Korean government(MSIT) (No. 2021R1A2C4002840).}

\maketitle

\section{Introduction}
\label{S1}

This paper deals with the Liouville-type theorem for the 3D stationary tropical climate model. 
The nonlinear partial differential equations
\begin{equation}
\label{E11}
\begin{split}
- \Delta u + (u \cdot \nabla ) u + \nabla \pi 
+ \divg (v \otimes v) &= 0, \\
- \Delta v + (u \cdot \nabla ) v + \nabla \theta 
+ (v \cdot \nabla )u &= 0, \\
- \Delta \theta + u \cdot \nabla \theta 
+ \divg v &= 0, \\
\divg u &= 0,
\end{split}
\end{equation}
in $\mathbb{R}^3$, describe the stationary tropical climate model.
Here, $u = (u_1,u_2,u_3)$ is the barotropic mode, $v = (v_1, v_2, v_3)$ is the first baroclinic mode of vector velocity, $\theta$ is the temperature, and $\pi$ is the pressure.

This paper aims to establish an improved Liouville-type theorem for the tropical climate model.
One of the most famous Liouville-type theorems is that if $f$ solves the Laplace equation on $\mathbb{R}^3$ and $f \in L^\infty(\mathbb{R}^3)$, then $f$ must be constant.
There are many variants of this theorem.
For example, if $f$ solves the Laplace equation on $\mathbb{R}^3$ and $f \in L^2(\mathbb{R}^3)$, then $f$ must be identically zero.
In general, Liouville-type theorems are about finding some conditions to show that solutions to some PDEs become trivial.
Recently, there have been many efforts to establish Liouville-type theorems for various fluid equations.
For the Navier–Stokes equations, one can find interesting results, for example, in  
Chae \cite{MR3162482}, Seregin \cite{MR3538409}, Kozono--Terasawa--Wakasugi \cite{MR3571910}, Chae--Wolf \cite{MR3959933}, and Cho--Choi--Yang \cite{MR4572397}.
For the tropical climate model, there are only a few results.
It was announced that a Liouville-type theorem holds if a smooth solution satisfies 
\begin{equation}
\label{E12}
u \in L^3 (\mathbb{R}^3), \quad 
v \in L^2 (\mathbb{R}^3), \quad \text{and} \quad
\nabla u, \nabla v, \nabla \theta \in L^2(\mathbb{R}^3),
\end{equation}
which is Theorem 1.3 in \cite{MR4342846}.
In the same paper, there are two other Liouville-type theorems.
We aim to remove the conditions $\nabla u, \nabla v \in L^2(\mathbb{R}^3)$.

Here is our main result.

\begin{thm}
\label{T1}
If a smooth solution $(u,v,\theta)$ to \eqref{E11} satisfies
\begin{equation}
\label{E13}
u \in L^3 (\mathbb{R}^3),
\quad
v \in L^2 (\mathbb{R}^3),
\quad \text{and} \quad
\nabla \theta \in L^2 (\mathbb{R}^3),
\end{equation}
then $u = v = 0$ and $\theta$ is constant.
\end{thm}

We derive an energy estimate, which provides a more useful direct proof of the Liouville-type property.
By using particular test functions with the Bogovskii operator and adapting an iteration method, we can remove the additional conditions $\nabla u, \nabla v \in L^2(\mathbb{R}^3)$ in \eqref{E12}. 
Indeed, we prove that the conditions \eqref{E13} imply $\nabla u, \nabla v \in L^2(\mathbb{R}^3)$.

\begin{rem}
Notice that one considers $\widetilde{\theta} = \theta + c$ for any constant $c$ instead of $\theta$ so that $\widetilde{\theta}$ solves the same PDEs and satisfies $\nabla \widetilde{\theta} \in L^2 (\mathbb{R}^3)$.
Hence, the conclusion that $\theta$ is constant in Theorem \ref{T1} is best possible.
\end{rem}

We end this section by giving a few notations and the Poincar\'e--Sobolev inequality used in this paper frequently.

\begin{itemize}
\item
For $0 < r < \infty$, we denote open balls and annuli by
\[
B(r) = \set{x \in \mathbb{R}^3 : 0 \le |x| < r}
\quad \text{and} \quad A(r) = \set{x \in \mathbb{R}^3 : r/2 < |x| < r}.
\]
\item
We will denote the Lebesgue measure of a measurable set $\Omega \subset \mathbb{R}^3$ by $|\Omega|$ and the Lebesgue integral of $f$ over $\Omega$ by $\int_\Omega f = \int_\Omega f(x) dx$.
\item
We will denote $L_0^p(\Omega)  = \set{f \in L^p (\Omega) : f_\Omega = 0}$, where the average value of $f$ over $\Omega$ is given by $f_\Omega = \frac{1}{|\Omega|} \int_\Omega f$.
\item
We will denote $A \lesssim B$ if $|A| \le c |B|$ for a generic positive constant $c$.
\end{itemize}

The following Lemma is called the Poincar\'e--Sobolev inequality.

\begin{lem}{[Theorem 3.15, \cite{MR1962933}]}
\label{L11}
Let $\Omega \subset \mathbb{R}^n$ be a bounded connected open set with Lipschitz-continuous boundary $\partial \Omega$.
There exists a positive constant $c (n,p,\Omega)$ such that if $p < n$, then we have for every $f \in W^{1,p} (\Omega)$,
\[
\norm{f - f_{\Omega}}_{L^\frac{np}{n-p} (\Omega)}
\le c (n,p,\Omega) \norm{\nabla f}_{L^p (\Omega)}.
\]
\end{lem}

\begin{rem}
[the Poincar\'e--Sobolev inequality on annuli]
If $\Omega = A(r)$ in the previous lemma, then the constant $c (n,p,\Omega)$ does not depend on $r>0$.
One can easily verify this by using a scaling method.
In particular, if we fix $n=3$ and $p=2$, then there is an absolute positive constant $c$ such that
\begin{equation}
\label{E14}
\norm{f - f_{A (r)}}_{L^6 (A (r))}
\le c \norm{\nabla f}_{L^2 ( A (r))}.
\end{equation}
\end{rem}

\section{Proof of Theorem \ref{T1}}
\label{S2}

We divide the proof into a few steps. 

\begin{enumerate}[Step 1.]
\item
(Derive an energy estimate) \\
We may assume that 
\begin{equation}
\label{E21}
\max\set{\norm{u}_3, \norm{v}_2, \norm{\nabla \theta}_2} \le M < \infty.
\end{equation}
Let $1 < R \le \rho < r \le 2R < \infty$ and $\varphi_{\rho,r} \in C_c^{\infty} (B(r))$ be a radially decreasing function such that $\varphi_{\rho,r} = 1$ on $B(\rho)$ and 
\begin{equation}
\label{E22}
(r-\rho) |\nabla \varphi_{\rho,r}| + (r-\rho)^2 |\nabla^2 \varphi_{\rho,r}| \le N < \infty,
\end{equation}
where $N$ is an absolute constant.
Using the Bogovskii operator $\mathcal{B}$, we can define
\[
w = \mathcal{B} (u \cdot \nabla \varphi_{\rho,r})
\]
in $A(r)$ since the support of $\nabla \varphi_{\rho,r}$ is contained in $A(r)$ and $u \cdot \nabla \varphi_{\rho,r} \in L_0^p(A(r))$.
Then $\divg w = u \cdot \nabla \varphi_{\rho,r}$ and for $1 < p < \infty$
\begin{equation}
\label{E23}
\norm{\nabla w}_{L^p (A(r))} \lesssim (r-\rho)^{-1} \norm{u}_{L^p (A(r))},
\end{equation}
where the implied constant depends only on $p$
(see \cite[Lemma 3]{MR4354995} or \cite[Lemma 4]{MR4572397} for the properties of the Bogovskii operator).
We will use the Einstein summation convention to sum over repeated indices.
We multiply the first equation of \eqref{E11} by $(u \varphi_{\rho,r} - w)$, the second equation of \eqref{E11} by $v \varphi_{\rho,r}$, and the third equation of \eqref{E11} by $(\theta - \theta_{A(r)}) \varphi_{\rho,r}$, and then integrate by parts with $\divg u = 0$ to obtain 
\begin{align*}
\int |\nabla u|^2 \varphi_{\rho,r}
&= \frac{1}{2} \int |u|^2 \Delta \varphi_{\rho,r} 
+ \frac{1}{2} \int |u|^2 u \cdot \nabla \varphi_{\rho,r} 
+ \int \partial_i u_j \partial_i w_j - \int u_i u_j \partial_i w_j \\
&\quad + \int v_i v_j \partial_i u_j \varphi_{\rho,r} 
+ \int v_i v_j \partial_i \varphi_{\rho,r} u_j - \int v_i v_j \partial_i w_j, \\
\int |\nabla v|^2 \varphi_{\rho,r}
&= \frac{1}{2} \int |v|^2 \Delta \varphi_{\rho,r} 
+ \frac{1}{2} \int |v|^2 u \cdot \nabla \varphi_{\rho,r} 
- \int v_i v_j \partial_i u_j \varphi_{\rho,r} \\
&\quad + \int \partial_j v_j (\theta - \theta_{A(r)}) \varphi_{\rho,r} 
+ \int (\theta - \theta_{A(r)}) v_j \partial_j \varphi_{\rho,r}, \\
\int |\nabla \theta|^2 \varphi_{\rho,r}
&= \frac{1}{2} \int |\theta - \theta_{A(r)}|^2 \Delta \varphi_{\rho,r} 
+ \frac{1}{2} \int |\theta - \theta_{A(r)}|^2 u \cdot \nabla \varphi_{\rho,r} 
- \int \partial_j v_j (\theta - \theta_{A(r)}) \varphi_{\rho,r}.
\end{align*}
Adding these identities, the four terms on the right are canceled so that we get 
\begin{align*}
&\int \left(|\nabla u|^2 + |\nabla v|^2 + |\nabla \theta|^2\right) \varphi_{\rho,r} \\
&= \frac{1}{2} \int \left(|u|^2 + |v|^2 + |\theta - \theta_{A(r)}|^2\right) \Delta \varphi_{\rho,r} 
+ \frac{1}{2} \int \left(|u|^2 + |v|^2 + |\theta - \theta_{A(r)}|^2\right) u \cdot \nabla \varphi_{\rho,r} 
+ \int v_i v_j \partial_i \varphi_{\rho,r} u_j \\
&\quad
+ \int (\theta - \theta_{A(r)}) v_j \partial_j \varphi_{\rho,r}
+ \int \partial_i u_j \partial_i w_j
- \int u_i u_j \partial_i w_j
- \int v_i v_j \partial_i w_j.
\end{align*}
Using \eqref{E22}, we have 
\begin{equation}
\label{E24}
\begin{split}
\int_{B(\rho)} \left( |\nabla u|^2 + |\nabla v|^2 + |\nabla \theta|^2 \right) 
&\lesssim (r-\rho)^{-2} \int_{A(r)} \left(|u|^2 + |v|^2 + |\theta - \theta_{A(r)}|^2\right) \\
&\quad + (r-\rho)^{-1} \int_{A(r)} \left(|u|^3 + |u||v|^2 + |v||\theta - \theta_{A(r)}| + |u||\theta - \theta_{A(r)}|^2\right) \\
&\quad + \int_{A(r)} |\nabla u| |\nabla w| 
+ \int_{A(r)} \left( |u|^2 |\nabla w| + |v|^2 |\nabla w| \right).
\end{split}
\end{equation}
\item
(Set up for an iteration) \\
Using the H\"older inequality, the Poincar\'e  inequality, and $1 < r \le 2R$, we get 
\begin{equation}
\label{E25}
\int_{A(r)} \left(|u|^2 + |v|^2 + |\theta - \theta_{A(r)}|^2\right) 
\lesssim r \norm{u}_3^2 + \norm{v}_2^2 + r^2 \norm{\nabla \theta}_{L^2 (A(r))}^2
\lesssim R + R^2 \norm{\nabla \theta}_{L^2 (A(r))}^2.
\end{equation}
Similarly, by the H\"older inequality, \eqref{E21}, and the Poincar\'e--Sobolev inequality \eqref{E14}, 
\begin{equation}
\label{E26}
\begin{split}
&\int_{A(r)} \left(|u|^3 + |u||v|^2 + |v||\theta - \theta_{A(r)}| + |u||\theta - \theta_{A(r)}|^2\right) \\
&\lesssim \norm{u}_3^3 
+ \norm{u}_3 \norm{v}_2 \norm{v}_{L^6(A(r))} 
+ r \norm{v}_2 \norm{\theta - \theta_{A(r)}}_{L^6(A(r))} 
+ r \norm{u}_3 \norm{\theta - \theta_{A(r)}}_{L^6(A(r))}^2 \\
&\lesssim 
1 + \norm{v}_{L^6(A(r))} 
+ R \norm{\nabla \theta}_{L^2 (A(r))}
+ R \norm{\nabla \theta}_{L^2 (A(r))}^2.
\end{split}
\end{equation}
By the H\"older inequality, \eqref{E23}, and \eqref{E21}, we obtain that
\begin{equation}
\label{E27}
\begin{split}
\int_{A(r)} |\nabla u| |\nabla w| 
&\le r^{1/2} \norm{\nabla u}_{L^2(A(r))} \norm{\nabla w}_{L^3(A(r))} 
\lesssim (r-\rho)^{-1} R^{1/2} \norm{\nabla u}_{L^2 (A(r))} \norm{u}_{L^3 (A(r))} \\
&\lesssim (r-\rho)^{-1} R^{1/2} \norm{\nabla u}_{L^2 (A(r))}
\end{split}
\end{equation}
and 
\begin{equation}
\label{E28}
\begin{split}
\int_{A(r)} \left( |u|^2 |\nabla w| + |v|^2 |\nabla w| \right) 
&\le \norm{u}_3^2 \norm{\nabla w}_{L^3(A(r))}
+ \norm{v}_2 \norm{v}_{L^6 (A(r))} \norm{\nabla w}_{L^3(A(r))} \\
&\lesssim (r-\rho)^{-1} \norm{u}_3^2 \norm{u}_{L^3 (A(r))}
+ (r-\rho)^{-1} \norm{v}_2 \norm{v}_{L^6 (A(r))} \norm{u}_{L^3 (A(r))} \\
&\lesssim (r-\rho)^{-1} + (r-\rho)^{-1} \norm{v}_{L^6 (A(r))}.
\end{split}
\end{equation}
By the Poincar\'e--Sobolev inequality \eqref{E14}, the Jensen inequality, and \eqref{E21}
\begin{equation}
\label{E29}
\norm{v}_{L^6 (A(r))}
\le \norm{v-v_{A(r)}}_{L^6 (A(r))}
+ \norm{v_{A(r)}}_{L^6 (A(r))} 
\lesssim \norm{\nabla v}_{L^2 (A(r))}
+ r^{-1} \norm{v}_{L^2 (A(r))}
\lesssim \norm{\nabla v}_{L^2 (A(r))} + 1. 
\end{equation}
Combining the estimates \eqref{E24}--\eqref{E29} gives
\begin{equation}
\label{E210}
\begin{split}
&\int_{B(\rho)} \left( |\nabla u|^2 + |\nabla v|^2 + |\nabla \theta|^2 \right) \\
&\lesssim (r-\rho)^{-2} (R + R^2 \norm{\nabla \theta}_{L^2 (A(r))}^2) \\
&\quad + (r-\rho)^{-1} (1 + \norm{\nabla v}_{L^2 (A(r))}
+ R \norm{\nabla \theta}_{L^2 (A(r))}
+ R \norm{\nabla \theta}_{L^2 (A(r))}^2) \\
&\quad + (r-\rho)^{-1} R^{1/2} \norm{\nabla u}_{L^2 (A(r))} 
+ (r-\rho)^{-1} + (r-\rho)^{-1} \norm{\nabla v}_{L^2 (A(r))}.
\end{split}
\end{equation}
Since $\nabla \theta \in L^2(\mathbb{R}^3)$, we have by the Young  inequality
\begin{align*}
\int_{B(\rho)} \left(|\nabla u|^2 + |\nabla v|^2\right) 
&\lesssim (r-\rho)^{-2} R^2 
+ (r-\rho)^{-1} \norm{\nabla v}_{L^2 (A(r))}
+ (r-\rho)^{-1} R \norm{\nabla u}_{L^2 (A(r))} \\
&\le \frac{1}{2} \int_{B(r)} \left(|\nabla u|^2 + |\nabla v|^2\right) 
+ c R^2 (r-\rho)^{-2}
\end{align*}
for some absolute constant $c>0$.
\item
(Vanishing energies at infinity) \\
We can apply the standard iteration argument (see \cite[Lemma 2]{MR4591749} or \cite[V. Lemma 3.1]{MR717034}) to obtain that for all $R \le \rho < r \le 2R$,
\[
\int_{B(\rho)} \left(|\nabla u|^2 + |\nabla v|^2\right) \le c R^2 (r-\rho)^{-2}.
\]
We now choose $\rho = R$ and $r = 2R$ so that 
\[
\int_{B(R)} \left(|\nabla u|^2 + |\nabla v|^2\right) \le c.
\]
Letting $R \to \infty$, we get $\nabla u, \nabla v \in L^2 (\mathbb{R}^3)$ and 
\begin{equation}
\label{E211}
\lim_{R \to \infty}
\int_{A(2R)} \left(|\nabla u|^2 + |\nabla v|^2 + |\nabla \theta|^2 \right) = 0.
\end{equation}
\item
(Vanishing energies on the whole space)\\
Using \eqref{E210} with $\rho = R$, $r = 2R$, we obtain
\begin{align*}
&\int_{B(R)} \left( |\nabla u|^2 + |\nabla v|^2 + |\nabla \theta|^2 \right) \\
&\lesssim R^{-2} (R + R^2 \norm{\nabla \theta}_{L^2 (A(2R))}^2) 
+ R^{-1} (1 + \norm{\nabla v}_{L^2 (A(2R))}
+ R \norm{\nabla \theta}_{L^2 (A(2R))}
+ R \norm{\nabla \theta}_{L^2 (A(2R))}^2) \\
&\quad + R^{-1/2} \norm{\nabla u}_{L^2 (A(2R))} 
+ R^{-1} + R^{-1} \norm{\nabla v}_{L^2 (A(2R))} \\
&\lesssim R^{-1} + \norm{\nabla \theta}_{L^2 (A(2R))}^2
+ R^{-1} \norm{\nabla v}_{L^2 (A(2R))}
+ \norm{\nabla \theta}_{L^2 (A(2R))}
+ R^{-1/2} \norm{\nabla u}_{L^2 (A(2R))}.
\end{align*}
Letting $R \to \infty$ and using \eqref{E211}, we conclude that 
\[
\lim_{R \to \infty}
\int_{B(R)} \left(|\nabla u|^2 + |\nabla v|^2 + |\nabla \theta|^2 \right) = 0,
\]
which gives $\nabla u = \nabla v = \nabla \theta = 0$.
Hence $u, v, \theta$ are constant.
Since $u \in L^3 (\mathbb{R}^3)$ and $v \in L^2 (\mathbb{R}^3)$, we should have $u = v = 0$. 
This completes the proof of Theorem \ref{T1}.
\end{enumerate}

\bibliography{my}

\bibliographystyle{plain}

\end{document}